\documentclass[12pt]{article}
\usepackage{amssymb,amsmath,mathrsfs}
\usepackage{hyperref}
\usepackage{graphicx}
\usepackage{color}
\usepackage[OT2,OT1]{fontenc}
\usepackage{upquote}
\usepackage[numbers,sort&compress]{natbib}

\begin{document}

\title{\LARGE\bf A rational approximation of the sinc function based on sampling and
\\ the Fourier transforms}

\author{
\normalsize\bf S. M. Abrarov\footnote{\scriptsize{Dept. Earth and Space Science and Engineering, York University, Toronto, Canada, M3J 1P3.}}\, and B. M. Quine$^{*}$\footnote{\scriptsize{Dept. Physics and Astronomy, York University, Toronto, Canada, M3J 1P3.}}}

\date{October 22, 2019}
\maketitle

\begin{abstract}
In our previous publications we have introduced the cosine product-to-sum identity \cite{Quine2013} 
$$
\prod\limits_{m = 1}^M {\cos \left( {\frac{t}{{{2^m}}}} \right)}  = \frac{1}{{{2^{M - 1}}}}\sum\limits_{m = 1}^{{2^{M - 1}}} {\cos \left( {\frac{{2m - 1}}{{{2^M}}}t} \right)}
$$
and applied it for sampling \cite{Abrarov2015a, Abrarov2015b} as an incomplete cosine expansion of the sinc function in order to obtain a rational approximation of the Voigt/complex error function that with only $16$ summation terms can provide accuracy ${\sim 10^{ - 14}}$. In this work we generalize this approach and show as an example how a rational approximation of the sinc function can be derived. A MATLAB code validating these results is presented.
\vspace{0.25cm}
\vspace{0.25cm}
\end{abstract}

\section{Introduction}

Sampling is a powerful mathematical tool that can be utilized in many fields of Applied Mathematics and Computational Physics \cite{Weisstein2003, Stenger2011, Rybicki1989}. One of the popular methods of sampling is based on the sinc function that can be defined as \cite{Kac1959, Gearhart1990}
$$
{\rm{sinc}}\left( t \right) = \left\{ 
\begin{aligned}
&\frac{{\sin t}}{t}, \quad\,\, t \ne 0\\
&1, \qquad\quad t = 0.
\end{aligned} \right.
$$
In particular, any function $f\left( t \right)$ defined in some interval $t \in \left[ {a,b} \right]$ can be approximated by the following equation (see, for example, equation (3) in \cite{Rybicki1989})
\begin{equation}\label{eq_1}
f\left( t \right) = \sum\limits_{n =  - N}^N {f\left( {{t_n}} \right){\text{sinc}}} \left( {\frac{\pi }{h}\left( {t - {t_n}} \right)} \right) + \varepsilon \left( t \right),
\end{equation}
where ${t_n}$ are the sampling points, $h$ is the small adjustable parameter and $\varepsilon \left( t \right)$ is the error term.

The forward and inverse Fourier transforms of function $f\left(t\right)$ can be defined as \cite{Bracewell2000, Hansen2014}
\begin{equation}\label{eq_2}
\mathcal{F}\left\{ {f\left( t \right)} \right\}\left( \nu  \right) = F\left( \nu  \right) = \int\limits_{ - \infty }^\infty  {f\left( t \right){e^{ - 2\pi i\nu t}}dt}
\end{equation}
and
\begin{equation}\label{eq_3}
{\mathcal{F}^{ - 1}}\left\{ {F\left( \nu  \right)} \right\}\left( t \right) = f\left( t \right) = \int\limits_{ - \infty }^\infty  {F\left( \nu  \right){e^{2\pi i\nu t}}d\nu },
\end{equation}
respectively. In our earlier works we have applied sampling by incomplete cosine expansion of the sinc function and obtained a rapid and highly accurate rational approximation for the Voigt/complex error function that with only $16$ summation terms can provide accuracy $\sim 10^{-14}$ \cite{Abrarov2015a, Abrarov2015b}. As a further development, in this work we generalize the methodology described in our paper \cite{Abrarov2015a} by an example showing how sampling and the Fourier transforms \eqref{eq_2}, \eqref{eq_3} can be implemented to obtain a rational approximation of the sinc function.

\section{Derivation}

The French mathematician Fran\c{c}ois Vi\`{e}te discovered an elegant formula relating the sinc function with cosine infinite product that can be written as given by \cite{Kac1959, Gearhart1990}
\begin{equation}\label{eq_4}
{\text{sinc}}\left( t \right) = \prod\limits_{m = 1}^\infty  {\cos \left( {\frac{t}{{{2^m}}}} \right)}.
\end{equation}

In our earlier publications we have applied the following cosine product-to-sum identity\footnote{Initially this identity was reported in our work \cite{Quine2013} to resolve a different problem in numerical integration related to the spectrally integrated Voigt function.}
\begin{equation}\label{eq_5}
\prod\limits_{m = 1}^M {\cos \left( {\frac{t}{{{2^m}}}} \right)}  = \frac{1}{{{2^{M - 1}}}}\sum\limits_{m = 1}^{{2^{M - 1}}} {\cos \left( {\frac{{2m - 1}}{{{2^M}}}t} \right)}
\end{equation}
for sampling to expand the Gaussian function $e^{-t^2}$ for numerical integration \cite{Abrarov2015a, Abrarov2015b}. As we can see, the left side of this identity represents a truncated cosine product \eqref{eq_4}. Consequently, the right side of the identity \eqref{eq_5} can be regarded as an incomplete cosine expansion of the sinc function. It is interesting to note that this identity has also found some useful applications in Computational Finance \cite{Ortiz-Gracia2016, Colldeforns-Papiol2018, Maree2018}. Here we show how it can also be applied to derive a rational approximation of the sinc function.

Due to product-to-sum identity \eqref{eq_5} the sinc function \eqref{eq_4} can be rewritten as the following limit
\[
{\text{sinc}}\left( t \right) = \mathop {\lim }\limits_{M \to \infty } \frac{1}{{{2^{M - 1}}}}\sum\limits_{m = 1}^{{2^{M - 1}}} {\cos \left( {\frac{{2m - 1}}{{{2^M}}}t} \right)}.
\]
Consequently, change of the variable $t \to \pi t/h$ in this limit yields
\[
{\text{sinc}}\left( {\frac{\pi }{h}t} \right) = \mathop {\lim }\limits_{M \to \infty } \frac{1}{{{2^{M - 1}}}}\sum\limits_{m = 1}^{{2^{M - 1}}} {\cos \left( {\frac{{\pi \left( {2m - 1} \right)}}{{{2^M}h}}t} \right)}.
\]
As we can see now, the truncated form of this limit can be implemented for sampling in accordance with equation \eqref{eq_1}. Specifically, the sinc function can be approximated as a cosine series expansion such that
\begin{equation}\label{eq_6}
{\text{sinc}}\left( {\frac{\pi }{h}t} \right)\approx \frac{1}{{{2^{M - 1}}}}\sum\limits_{m = 1}^{{2^{M - 1}}} {\cos \left( {\frac{{\pi \left( {2m - 1} \right)}}{{{2^M}h}}t} \right)}, \qquad - T/4 \leqslant t \leqslant T/4,
\end{equation}
where $T = {2^{M + 1}}h$. Since equation \eqref{eq_6} consists of fixed number of the cosine terms, it is, therefore, periodic as it is shown in Fig. 1.

\begin{figure}[ht]
\begin{center}
\includegraphics[width=24pc]{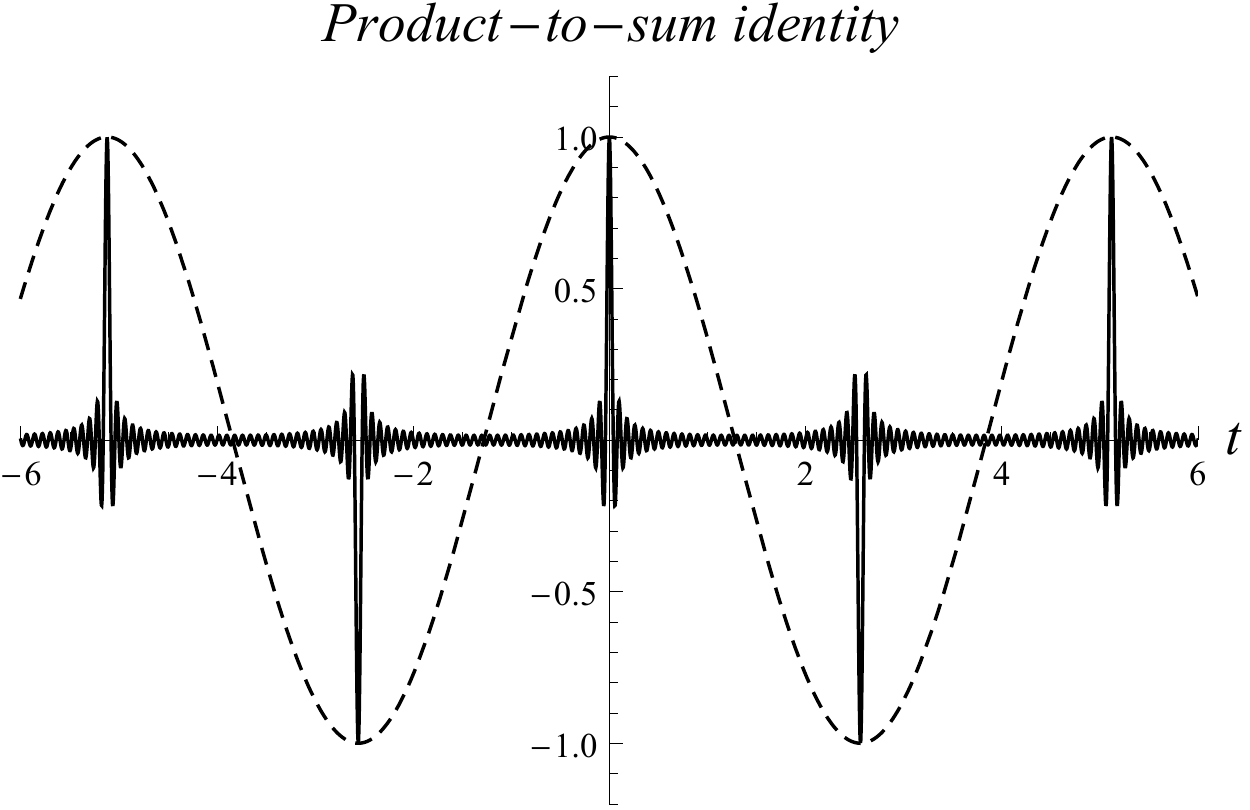}\hspace{2pc}%
\begin{minipage}[b]{28pc}
\vspace{0.3cm}
{\sffamily {\bf{Fig. 1.}} The product-to-sum identity \eqref{eq_5} upon change of the variable $t\to \pi t/h$ at $M=6$ and $h=0.04$. The dashed curve shows ${\cos \left(\frac{\pi t}{2^6 \times 0.04} \right)}$ corresponding to the first summation term of approximation \eqref{eq_6}.}
\end{minipage}
\end{center}
\end{figure}

Generally, the sampling grid-points may not be equidistant. However, for simplicity we can assume equidistantly spaced grid-points along $t$-axis. Thus, by taking ${t_n} = nh$ the substitution of approximation \eqref{eq_6} into equation \eqref{eq_1} results in
\footnotesize
\begin{equation}\label{eq_7}
f\left( t \right) \approx \frac{1}{{{2^{M - 1}}}}\sum\limits_{m = 1}^{{2^{M - 1}}} {\sum\limits_{n =  - N}^N {f\left( {nh} \right)\cos \left( {\frac{{\pi \left( {2m - 1} \right)}}{{{2^M}h}}\left( {t - nh} \right)} \right)} } , \quad - T/4 \leqslant t \leqslant T/4.
\end{equation}
\normalsize
Similar to equation \eqref{eq_6} the right side of equation \eqref{eq_7} is also periodic and, consequently, it approximates only within the range $t \in \left[ { - T/4,T/4} \right]$.

At the beginning we need to determine the inverse Fourier transform for the sinc function ${\text{sinc}}\left(\pi \nu \right)$ by using equation \eqref{eq_3}. This leads to
$$
{\mathcal{F}^{ - 1}}\left\{ {{\text{sinc}}\left( {\pi \nu } \right)} \right\}\left( t \right) = \int\limits_{ - \infty }^\infty  {{\text{sinc}}\left( {\pi \nu } \right){e^{2\pi i\nu t}}d\nu }  = {\text{rect}}\left( t \right),
$$
where
\[
{\text{rect}}\left( t \right) = \left\{
\begin{aligned} & 1, \qquad {\text{if}} \, \left|t\right| < 1/2
\\ & 1/2, \quad {\text{if}} \, \left|t\right| =  1/2
\\ & 0, \qquad {\text{if}} \, \left|t\right| > 1/2,
\end{aligned}  \right.
\]
is commonly known as the rectangular function.

The rectangular function contains two discontinuities at $t = - 1/2$ and $t = 1/2$. Therefore, it is somehow problematic to handle this function due to the high level of oscillation near these discontinuities that occurs as a result of sampling. However, from the well-known relation given by
\[
{\text{rect}}\left( t \right) = \mathop {\lim }\limits_{k \to \infty } \frac{1}{{{{\left( {2t} \right)}^{2k}} + 1}},
\]
where $k$ is a positive integer, it follows that the rectangular function ${\text{rect}}\left( t \right)$ can be approximated with high accuracy by taking sufficiently large integer $k$. Thus, if we take, say $k = 35$, then the rectangular function can be approximated very accurately as
\[
{\text{rect}}\left( t \right) \approx \frac{1}{{{{\left( {2t} \right)}^{70}} + 1}}.
\]
Although this approximation changes very rapidly near the points $t =  - 1/2$ and $t = 1/2$, it has no discontinuities and, therefore, can be used for sampling more efficiently then the function ${\text{rect}}\left( t \right)$ itself. Figure 2 shows the function $1/\left[ {{{\left( {2t} \right)}^{70}} + 1} \right]$ by blue curve.

\begin{figure}[ht]
\begin{center}
\includegraphics[width=24pc]{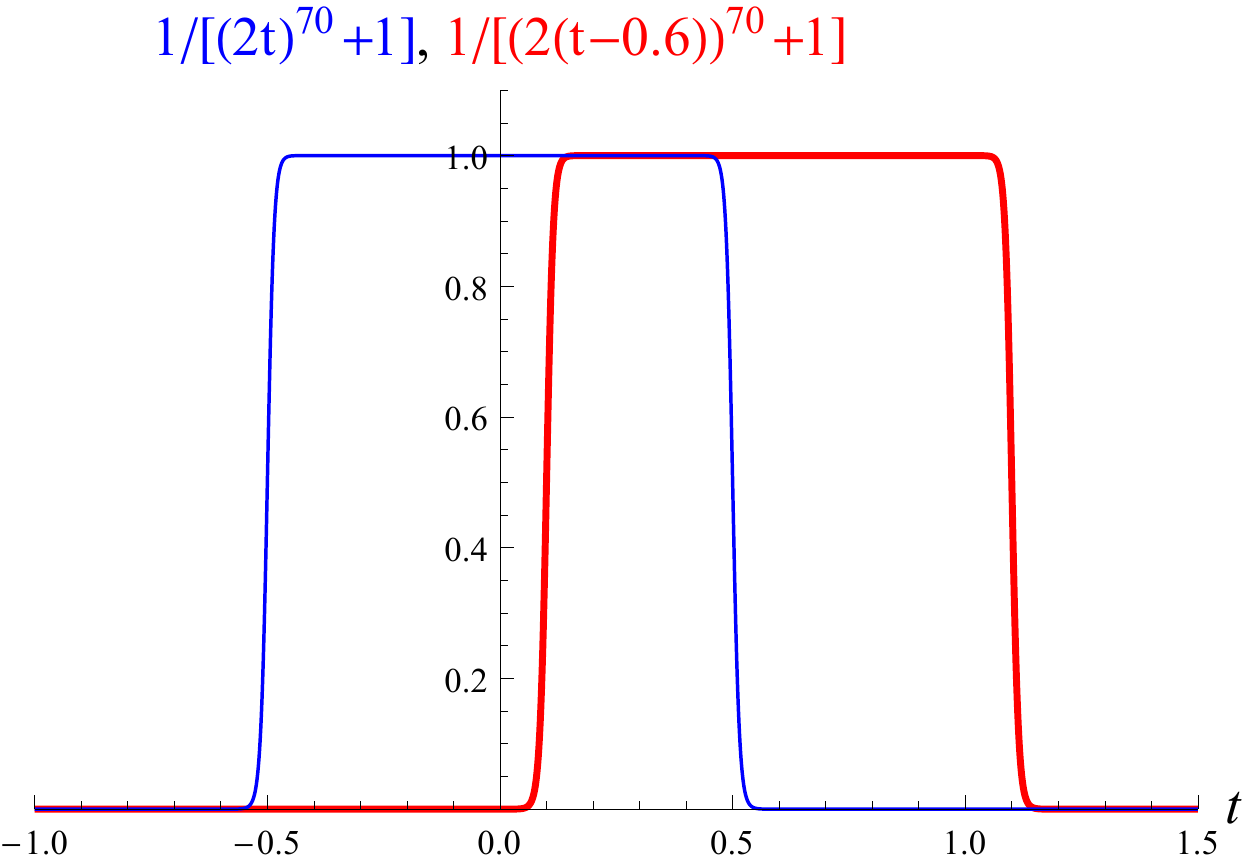}\hspace{2pc}%
\begin{minipage}[b]{28pc}
\vspace{0.3cm}
{\sffamily {\bf{Fig. 2.}} The rectangular function approximations $f\left( t \right) = 1/\left[ {{{\left( {2t} \right)}^{70}} + 1} \right]$ and $f\left( t - a \right) = 1/\left[ {{{\left( {2\left(t - 0.6\right)} \right)}^{70}} + 1} \right]$ shown by blue and red curves, respectively.}
\end{minipage}
\end{center}
\end{figure}

The function $1/\left[ {{{\left( {2t} \right)}^{70}} + 1} \right]$ effectively covers the range
$$
t \in \left[ { - \left( {1/2 + \Delta } \right),\left( {1/2 + \Delta } \right)} \right],
$$
where a small positive value $\Delta$ can be taken equal to $0.1$. By stating ``effectively covers the range'' we imply that the Fourier transform can be approximated as
$$
\mathcal{F}\left\{ \frac{1}{{{{\left( {2t} \right)}^{70}} + 1}} \right\}\left( \nu  \right) = \int\limits_{ - \infty }^\infty  {\frac{1}{{{{\left( {2t} \right)}^{70}} + 1}}{e^{ - 2\pi i\nu t}}dt}  \approx \int\limits_{ - a}^a {\frac{1}{{{{\left( {2t} \right)}^{70}} + 1}}{e^{ - 2\pi i\nu t}}dt},
$$
where $a = 1/2 + \Delta  = 0.6$, since
$$
\int\limits_{ - \infty }^{-a} {\frac{1}{{{{\left( {2t} \right)}^{70}} + 1}}{e^{ - 2\pi i\nu t}}dt}  \approx 0
$$
and
$$
\int\limits_a^\infty  {\frac{1}{{{{\left( {2t} \right)}^{70}} + 1}}{e^{ - 2\pi i\nu t}}dt}  \approx 0
$$
are negligibly small and, consequently, can be ignored.

In this approach we consider only the first quadrant. Therefore, we shift the function $f\left( t \right) = 1/\left[ {{{\left( {2t} \right)}^{70}} + 1} \right]$ towards right in form
\[
f\left( {t - a} \right) = \frac{1}{{{{\left( {2\left( {t - a} \right)} \right)}^{70}} + 1}} = \frac{1}{{{{\left( {2\left( {t - 0.6} \right)} \right)}^{70}} + 1}}
\]
as it is shown in Fig. 2 by red curve.

The Fourier transform approximated as
\begin{equation}\label{eq_8}
\mathcal{F}\left\{ {f\left( {t - a} \right)} \right\}\left( \nu  \right) = \int\limits_{ - \infty }^\infty  {f\left( {t - a} \right){e^{ - 2\pi i\nu t}}dt}  \approx \int\limits_0^{2a} {f\left( {t - a} \right){e^{ - 2\pi i\nu t}}dt} ,
\end{equation}
does not provide a rational approximation as a result of the fixed value $2a$ on the upper limit of the integration. However, this problem can be effectively resolved by sampling the function $f\left( {t - a} \right){e^{\sigma t}}$ instead of $f\left( {t - a} \right)$ as follows\footnote{If applicable, the Poisson summation formula can also be used instead of sampling to expand a function of kind $f\left( {t - a} \right){e^{\sigma t}}$. However, the sampling method is more advantageous due to its versatility in practical applications.}
\small
\begin{equation}\label{eq_9}
\begin{aligned}
f\left( {t - a} \right){e^{\sigma t}} \approx \frac{1}{{{2^{M - 1}}}}\sum\limits_{m = 1}^{{2^{M - 1}}} &{\sum\limits_{n =  - N}^N {f\left( {nh - a} \right){e^{\sigma nh}}\cos \left( {\frac{{\pi \left( {2m - 1} \right)}}{{{2^M}h}}\left( {t - nh} \right)} \right)} },
\\ & - T/4 + a \leqslant t \leqslant T/4 + a,
\end{aligned}
\end{equation}
\normalsize
where $\sigma $ is a positive (adjustable) number. Figure 3 illustrates the approximation \eqref{eq_9} to the function $f\left( t - a \right) e^{\sigma t} = e^{\sigma t}/\left[ {{{\left( {2\left(t-0.6\right)} \right)}^{70}} + 1} \right]$ at $a=0.6$, $M = 6$, $N = 28$, $h=0.04$ with $\sigma = 0$, $\sigma = 0.25$ and $\sigma = 0.75$ by blue, red and green curves, respectively.

\begin{figure}[ht]
\begin{center}
\includegraphics[width=26pc]{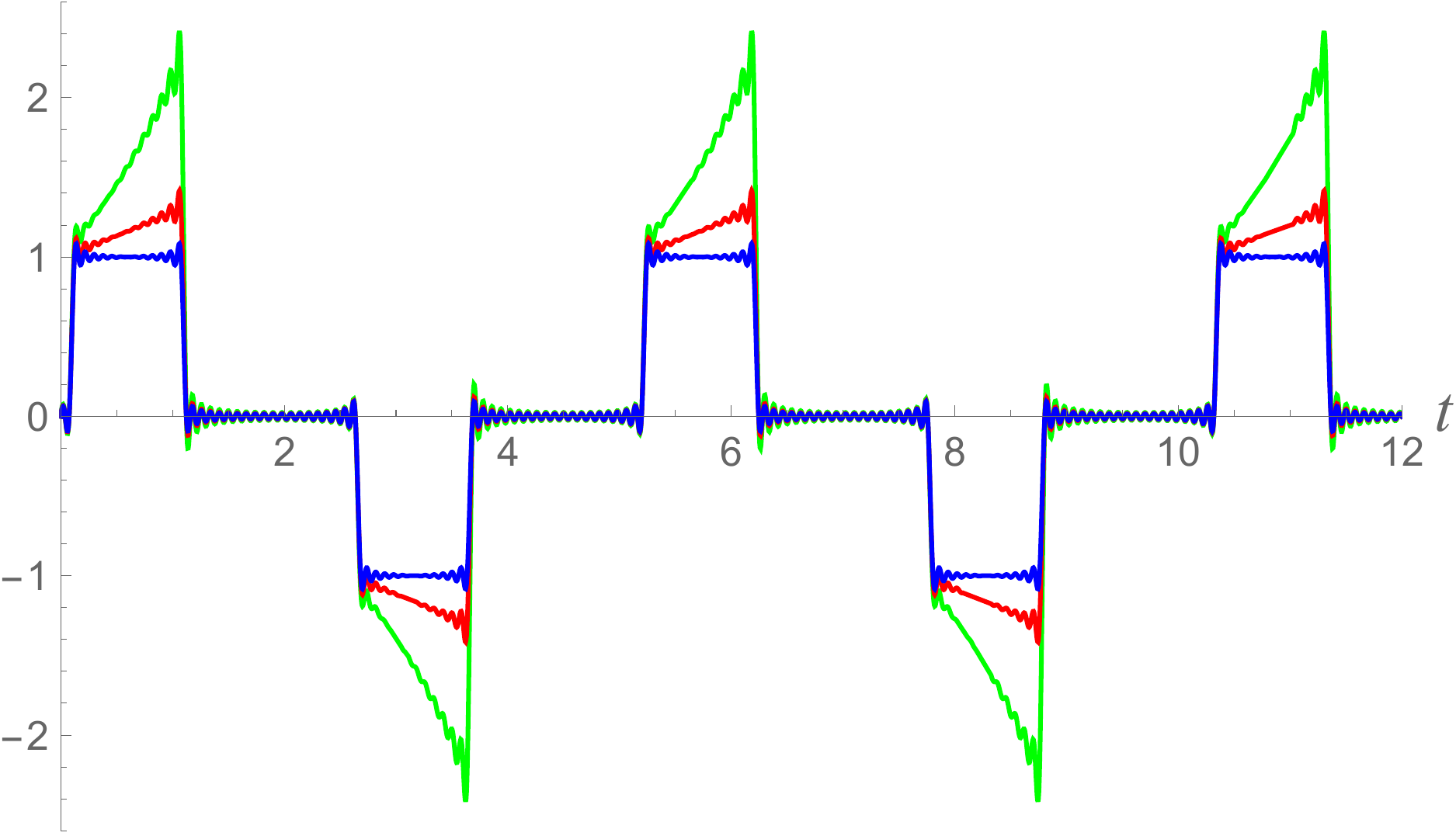}\hspace{2pc}%
\begin{minipage}[b]{28pc}
\vspace{0.3cm}
{\sffamily {\bf{Fig. 3.}} Approximation \eqref{eq_9} to the function $f\left( t - a \right) e^{\sigma t} = e^{\sigma t}/\left[ {{{\left( {2\left(t-0.6\right)} \right)}^{70}} + 1} \right]$ calculated at $a=0.6$, $M = 6$, $N = 28$, $h=0.04$ with $\sigma = 0$ (blue curve), $\sigma = 0.25$ (red curve) and $\sigma = 0.75$ (green curve).}
\end{minipage}
\end{center}
\end{figure}

Since the approximation \eqref{eq_9} is periodic, it remains valid only within the range $t \in \left[ { - T/4 + a,T/4 + a} \right]$. However, when $\sigma $ is big enough, say equal to or greater than $1$, the slight rearrangement of the approximation \eqref{eq_9} as given by
\footnotesize
\begin{equation}\label{eq_10}
f\left( {t - a} \right) \approx \frac{{{e^{ - \sigma t}}}}{{{2^{M - 1}}}}\sum\limits_{m = 1}^{{2^{M - 1}}} {\sum\limits_{n =  - N}^N {f\left( {nh - a} \right){e^{\sigma nh}}\cos \left( {\frac{{\pi \left( {2m - 1} \right)}}{{{2^M}h}}\left( {t - nh} \right)} \right)} } , \quad \sigma \gtrsim 1,
\end{equation}
\normalsize
extends the validation of the function $f\left( {t - a} \right)$ along entire positive $t$-axis. This is possible to achieve since the damping multiplier ${e^{ - \sigma t}}$ effectively suppresses all peaks (except the first peak near the origin) that appear as a result of this periodicity. This suppression effect given by the damping multiplier $e^{-\sigma t}$ can be seen from Fig. 4 that shows the evolution to the function $f\left( t - a \right)$ with increasing $\sigma$. In particular, at $\sigma  = 0$ we observe no suppression and the signal remains periodic (blue curve). However, as $\sigma $ increases the suppression effect due to presence of the damping multiplier ${e^{ - \sigma t}}$ increases with increasing $t$ as shown by red and green curves. Consequently, if the constant $\sigma$ is large enough, then the function becomes practically solitary. As a consequence, we can substitute approximation \eqref{eq_10} into the equation \eqref{eq_8} and replace the upper limit $2a$ in the integral by $\infty$ practically without loss of accuracy as follows
\footnotesize
\begin{equation}\label{eq_11}
\begin{aligned}
& \hspace{-0.5 cm} \mathcal{F}\left\{ {f\left( {t - a} \right)} \right\}\left( \nu  \right) 
\\ & \approx \int\limits_0^{2a} {\left[ {\frac{e^{-\sigma t}}{{{2^{M - 1}}}}\sum\limits_{m = 1}^{{2^{M - 1}}} {\sum\limits_{n =  - N}^N {f\left( {nh - a} \right){e^{\sigma nh}}\cos \left( {\frac{{\pi \left( {2m - 1} \right)}}{{{2^M}h}}\left( {t - nh} \right)} \right)} } } \right]{e^{ - 2\pi i\nu t}}dt}
\\ & \approx \int\limits_0^\infty  {\left[ {\frac{e^{-\sigma t}}{{{2^{M - 1}}}}\sum\limits_{m = 1}^{{2^{M - 1}}} {\sum\limits_{n =  - N}^N {f\left( {nh - a} \right){e^{\sigma nh}}\cos \left( {\frac{{\pi \left( {2m - 1} \right)}}{{{2^M}h}}\left( {t - nh} \right)} \right)} } } \right]{e^{ - 2\pi i\nu t}}dt}.
\end{aligned}
\end{equation}
\normalsize

Since the upper limit is infinity now the resultant integration yields a rational function. Thus, after some trivial rearrangements that exclude the double summation, from equation \eqref{eq_11} it follows that
\begin{equation}\label{eq_12}
\mathcal{F}\left\{ {f\left( {t - a} \right)} \right\}\left( \nu  \right) \approx \sum\limits_{m = 1}^{{2^{M - 1}}} {\frac{{{\alpha _m}\left( {\sigma + 2\pi i\nu} \right) + {\beta _m}}}{{\gamma _m^2 + {{\left( {\sigma + 2\pi i\nu} \right)}^2}}}},
\end{equation}
where the corresponding expansion coefficients are
\[
\label{eq_13a}
\tag{13a}
{\alpha _m} = \frac{1}{{{2^{M - 1}}}}\sum\limits_{n =  - N}^N {{f\left( {nh - a} \right){e^{\sigma nh}}\cos \left( {\frac{{\pi \left( {2m - 1} \right)}}{{{2^M}h}}nh} \right)} },
\]
\[
\label{eq_13b}
\tag{13b}
{\beta _m} = \frac{1}{{{2^{M - 1}}}}\sum\limits_{n =  - N}^N {{f\left( {nh - a} \right){e^{\sigma nh}}\frac{{\pi \left( {2m - 1} \right)}}{{{2^M}h}}\sin \left( {\frac{{\pi \left( {2m - 1} \right)}}{{{2^M}h}}nh} \right)}}
\]
and
\[
\label{eq_13c}
\tag{13c}
{\gamma _m} = \frac{{\pi \left( {2m - 1} \right)}}{{{2^M}h}}.
\]

\begin{figure}[ht]
\begin{center}
\includegraphics[width=24pc]{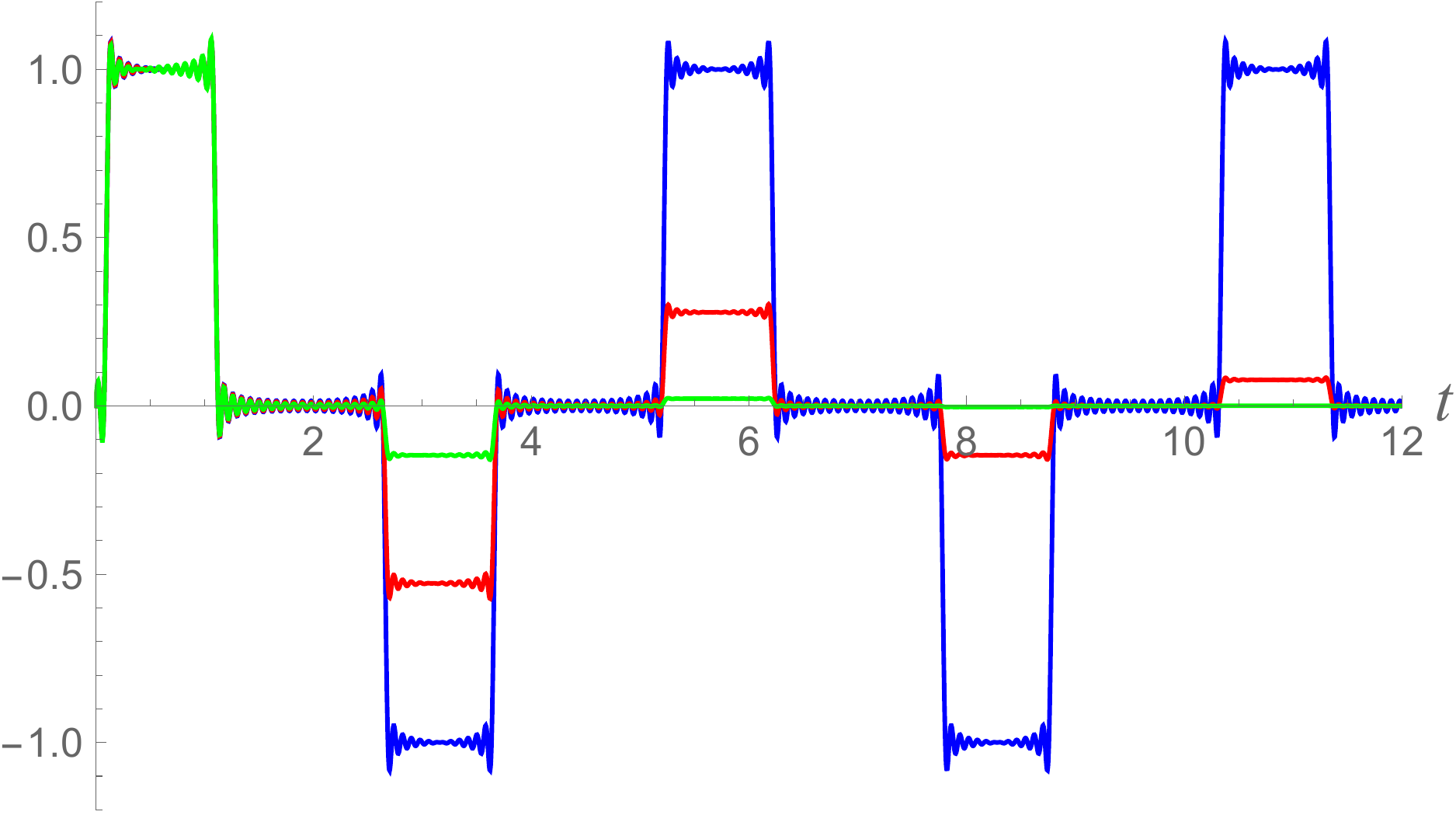}\hspace{2pc}%
\begin{minipage}[b]{28pc}
\vspace{0.3cm}
{\sffamily {\bf{Fig. 4.}} Evolution to the function $f\left( t - a \right) = 1/\left[ {{{\left( {2\left(t-0.6\right)} \right)}^{70}} + 1} \right]$ obtained from approximation \eqref{eq_10} at $a=0.6$, $M = 6$, $N = 28$, $h=0.04$ with $\sigma = 0$ (blue curve), $\sigma = 0.25$ (red curve) and $\sigma = 0.75$ (green curve).}
\end{minipage}
\end{center}
\end{figure}

As we can see, equations \eqref{eq_12}, \eqref{eq_13a}, \eqref{eq_13b} and \eqref{eq_13c} are of the same kind as that of reported in our earlier work \cite{Abrarov2015a}, where the damping multiplier $e^{-\sigma t}$ has also been used in numerical integration for derivation of the rational approximation of the Voigt/complex error function (compare Fig. 4 above with Fig. 3 from our paper \cite{Abrarov2015a} to observe same effect of suppression along positive $t$-axis due damping multiplier $e^{-\sigma t}$ ).

\section{Approximation}

Since the function $f\left( {t - a} \right)$ effectively covers the range from $0$ to $2a$ (see red curve in Fig. 2), it is not necessary to run the index $n$ starting from $ - N$. Therefore, the expansions \eqref{eq_13a} and \eqref{eq_13b} can be slightly modified by running the index $n$ from $0$ to $N$. This leads to
\[
\tag{14a}
\begin{aligned}{\alpha _m} &= \frac{1}{{{2^{M - 1}}}}\sum\limits_{n = 0}^N { {f\left( {nh - a} \right){e^{\sigma nh}}\cos \left( {\frac{{\pi \left( {2m - 1} \right)}}{{{2^M}h}}nh} \right)} }
\\ &= \frac{1}{{{2^{M - 1}}}}\sum\limits_{n = 0}^N { {f\left( {nh - a} \right){e^{\sigma nh}}\cos \left( {{\gamma _m}n h} \right)} } 
\end{aligned}
\]
and
\[
\tag{14b}
\begin{aligned}
{\beta _m} &= \frac{1}{{{2^{M - 1}}}}\sum\limits_{n = 0}^N {{f\left( {nh - a} \right){e^{\sigma nh}}\frac{{\pi \left( {2m - 1} \right)}}{{{2^M}h}}\sin \left( {\frac{{\pi \left( {2m - 1} \right)}}{{{2^M}h}}nh} \right)} }
\\ &= \frac{1}{{{2^{M - 1}}}}\sum\limits_{n = 0}^N {{f\left( {nh - a} \right){e^{\sigma nh}}{\gamma _m}\sin \left( {{\gamma _m}n h} \right)} }.
\end{aligned}
\]

Using the shifting property of the Fourier transform \eqref{eq_2} we can write
\addtocounter{equation}{2}
\begin{equation}\label{eq_15}
\mathcal{F}\left\{ {f\left( t \right)} \right\}\left( \nu  \right) = e^{2\pi i\nu a}\mathcal{F}\left\{ {f\left( {t - a} \right)} \right\}\left( \nu  \right).
\end{equation}
Consequently, from equations \eqref{eq_12} and \eqref{eq_15} we obtain
\begin{equation}\label{eq_16}
\begin{aligned}
\mathcal{F}\left\{ f \left( t \right) \right\} \left( \nu \right) &= \mathcal{F}\left\{ \frac{1}{{{{\left( {2t} \right)}^{70}} + 1}} \right\}\left( \nu  \right)
\\ &\approx e^{2\pi i\nu a}\sum\limits_{m = 1}^{{2^{M - 1}}} {\frac{{{\alpha _m} \left( {\sigma + 2\pi i\nu} \right) + {\beta _m}}}{{\gamma _m^2 + {{\left({\sigma + 2\pi i\nu} \right)}^2}}}} \approx \text{sinc} \left( \pi \nu \right).
\end{aligned}
\end{equation}

\begin{figure}[ht]
\begin{center}
\includegraphics[width=24pc]{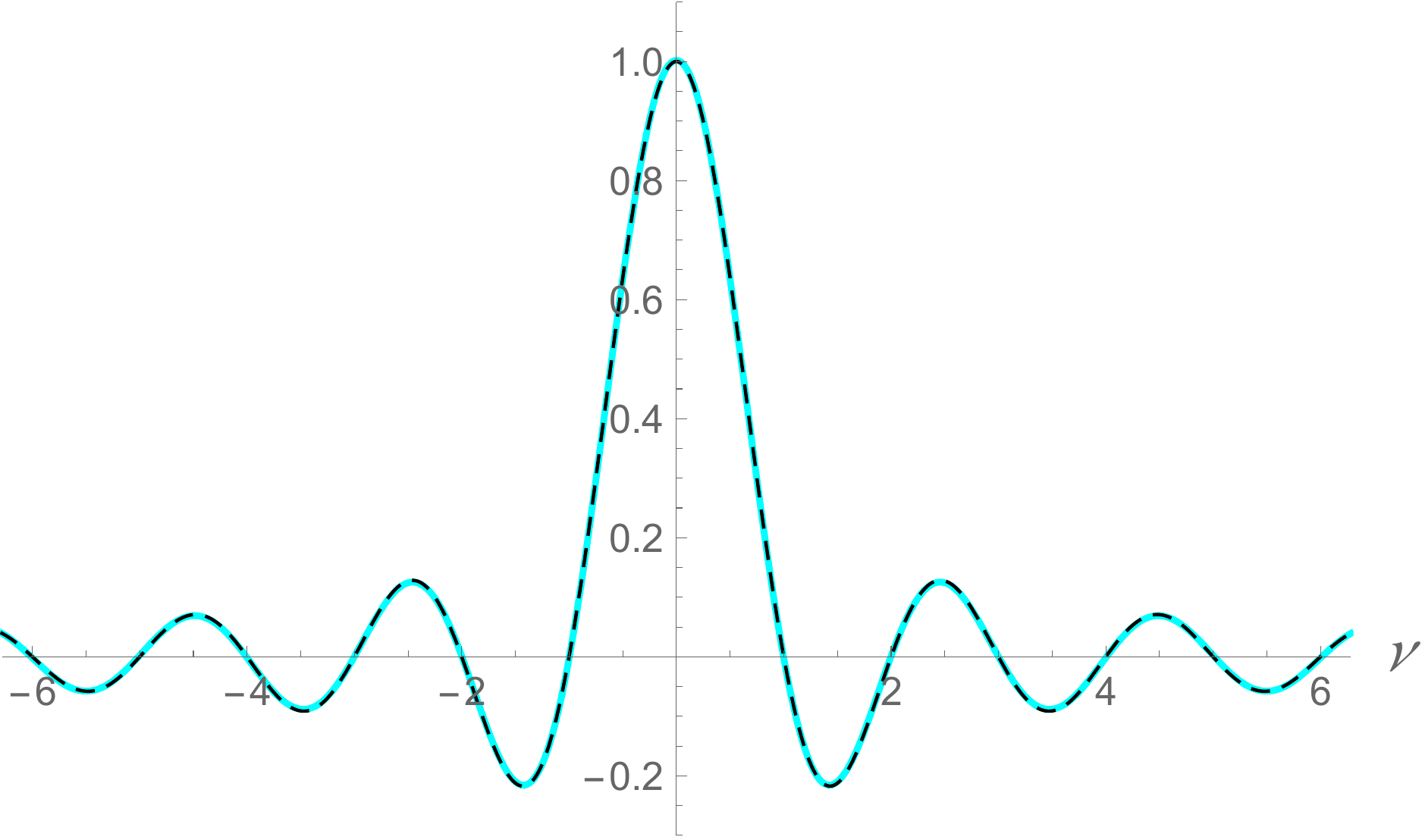}\hspace{2pc}%
\begin{minipage}[b]{28pc}
\vspace{0.3cm}
{\sffamily {\bf{Fig. 5.}} Sinc function approximation \eqref{eq_16} calculated at $a=0.6$, $M = 6$, $N = 28$, $h=0.04$ and $\sigma = 2.7$ (light blue curve). The original sinc function $\text{sinc}\left(\pi \nu \right)$ is also shown by dashed curve.}
\end{minipage}
\end{center}
\end{figure}

Figure 5 shows a rational approximation \eqref{eq_16} of the sinc function by light blue curve. The black dashed curve shows the original sinc function ${\text{sinc}}\left( {\pi \nu } \right)$ for comparison. As we observe, the original sinc function ${\text{sinc}}\left( {\pi \nu } \right)$ and its rational approximation \eqref{eq_16} are not distinctive visually.

\begin{figure}[ht]
\begin{center}
\includegraphics[width=24pc]{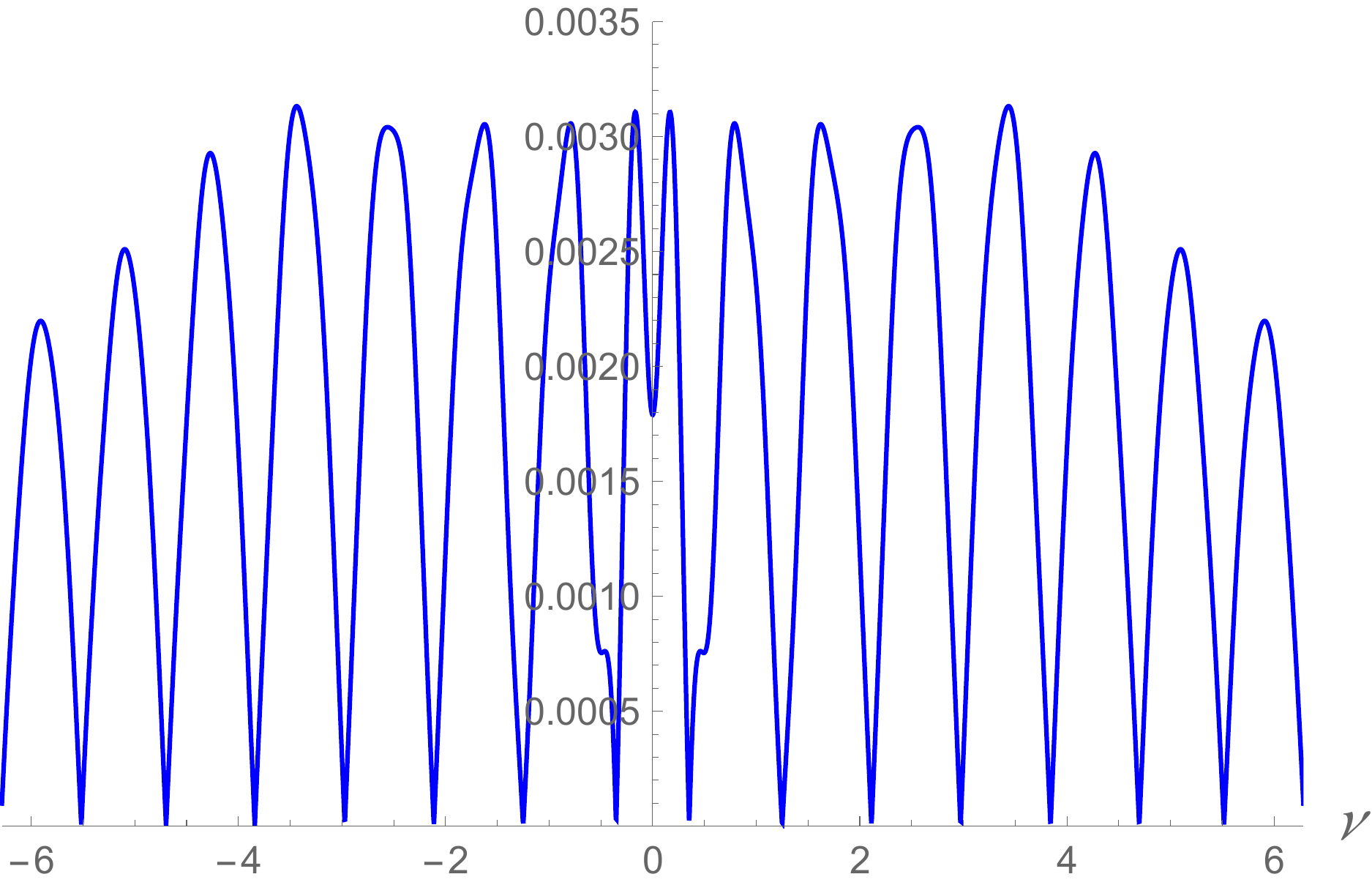}\hspace{2pc}%
\begin{minipage}[b]{28pc}
\vspace{0.3cm}
{\sffamily {\bf{Fig. 6.}} Absolute difference between the original sinc function $\text{sinc} \left( \pi \nu \right)$ and its approximation \eqref{eq_16} at $a=0.6$, $M = 6$, $N = 28$, $h=0.04$ and $\sigma = 2.7$.}
\end{minipage}
\end{center}
\end{figure}

Figure 6 shows the absolute difference between the original sinc function ${\text{sinc}}\left( {\pi \nu } \right)$ and its rational approximation. Despite that the sinc function is not easy to approximate \cite{Qiu2010, Lu2012}, only ${2^{6 - 1}} = 32$ terms in the proposed formula \eqref{eq_16} provide absolute difference smaller than $3.2 \times 10^{-3}$ within the range $\nu  \in \left[ { - 2\pi ,2\pi } \right]$. A MATLAB code validating the results based on the rational approximation \eqref{eq_16} of the sinc function ${\text{sinc}}\left( {\pi \nu } \right)$ is provided in Appendix A. Specifically, a command line {\sffamily{rationapp(1)}} displays two graphs with approximated sinc function and its absolute difference.

It should be noted that significantly higher accuracy can be achieved with more well-behaved functions. For example, substituting the function $f\left( t \right) = \pi^{3/2} i\,t e^{-\pi^2 t^2}$ into equation \eqref{eq_16} and taking parameters $a=2$, $M=6$, $N=55$, $h=0.078$ and $\sigma=5$ we can obtain a rational approximation of the function $\mathcal{F}\left\{ \pi^{3/2} i\,t e^{-\pi^2 t^2} \right\}\left( \nu \right) = \nu e^{-\nu ^2}$ with absolute difference less than $7.3 \times 10^{-12}$. This can be readily confirmed by running a MATLAB code provided in Appendix A. In particular, a sample computation that can be run by a command line {\sffamily{rationapp(2)}} shows that the absolute difference between original function $\nu e^{-\nu^2}$ and its rational approximation does not exceed $7 \times 10^{-12}$.

\section{Numerical integration}

The application of the rational approximation \eqref{eq_16} may be especially advantageous for a contour integral of kind
\begin{equation}\label{eq_17}
\begin{aligned}
&\hspace{-1cm}\oint\limits_C {\mathcal{F}\left\{ {f\left( t \right)} \right\}\left( \nu  \right) g\left( \nu, x,y,z,\dots \right) d\nu } \\
& \hspace{0.5cm} = \oint\limits_C { \left( e^{2\pi i\nu a}\sum\limits_{m = 1}^{{2^{M - 1}}} {\frac{{{\alpha _m} \left(\sigma + {2\pi i\nu } \right) + {\beta _m}}}{{\gamma _m^2 + {{\left( \sigma + {2\pi i\nu} \right)}^2}}}} \right)} g\left( \nu, x,y,z,\dots \right)\,d\nu,
\end{aligned}
\end{equation}
where $g\left( \nu, x,y,z,\dots \right)$ is a multivariable function, since it may be computed numerically by residues when the integrand on the right side of integral in equation \eqref{eq_17} is analytic. This is possible to achieve due to presence of the poles in the rational approximation of function $\mathcal{F}\left\{ {f\left( t \right)} \right\}\left( \nu  \right)$. For example, applying the functions
$$
f\left(t\right) = \sqrt{\pi} e^{-\pi^2 t^2} \Leftrightarrow \mathcal{F}\left\{\sqrt{\pi} e^{-\pi^2 t^2} \right\} \left( \nu  \right) = e^{-\nu^2}
$$
and
$$
g \left( \nu, x,y \right) = \frac{y}{{\pi\left(y^2+\left(x-\nu \right)^2\right)}}
$$
in equation \eqref{eq_17}, the following integral
\begin{equation}\label{eq_18}
K\left(x,y\right) =\frac{y}{\pi} \int\limits_{-\infty}^\infty \frac{e^{-\tau^2}}{y^2+\left(x-\tau \right)^2}  \,d\tau,
\end{equation}
known as the Voigt function \cite{Armstrong1967}, can be calculated by residues in a similar way that we performed in our earlier work \cite{Abrarov2015b} to approximate the integral \eqref{eq_18}. Due to rapid performance and high accuracy our algorithm \cite{Abrarov2015b} has been implemented in current version of the atmospheric model {\it{bytran}} \cite{Pliutau2017, bytran}. Thus, integration by residues leads to
\footnotesize
\begin{equation}\label{eq_19}
\begin{aligned}
  K\left( {x,y} \right) \approx 2\pi iy \sum\limits_{m = 1}^{{2^{M - 1}}} &{\left[{\frac{{{e^{ - a\left( {i{\gamma _m} + \sigma } \right)}}\left( {{\beta _m} - i{\alpha _m}{\gamma _m}} \right)}}{{{\gamma _m}\left( {4{\pi ^2}\left( {{x^2} + {y^2}} \right) + 4\pi x\left( {{\gamma _m} - i\sigma } \right) + {{\left( {{\gamma _m} - i\sigma } \right)}^2}} \right)}}} \right.}  \\ 
   &- \frac{{i{e^{a\left( {i{\gamma _m} - \sigma } \right)}}\left( {{\alpha _m}{\gamma _m} - i{\beta _m}} \right)}}{{{\gamma _m}\left( {4{\pi ^2}\left( {{x^2} + {y^2}} \right) - 4\pi x\left( {{\gamma _m} + i\sigma } \right) + {{\left( {{\gamma _m} + i\sigma } \right)}^2}} \right)}} \\ 
  &\left. { + \frac{{i{e^{2ia\pi \left( {x + iy} \right)}}\left( {{\alpha _m}\left( {2\pi \left( {y - ix} \right) - \sigma } \right) - {\beta _m}} \right)}}{{2\pi y\left( {\gamma _m^2 - {{\left( {2\pi \left( {x + iy} \right) - i\sigma } \right)}^2}} \right)}}} \right]. 
\end{aligned} 
\end{equation}
\normalsize
This approximation can provide highly accurate values as we can see from the Fig. 7 showing an example of the absolute difference between the original Voigt function $K\left(x,y\right)$ and its approximation \eqref{eq_19}. In particular, a computation performed with extended precision floating-point by using Wolfram Mathematica (version 11) shows that at $y=1$, $a=2$, $M=6$, $N=55$, $h=0.078$ and $\sigma = 5$ the absolute difference does not exceed $4 \times 10^{-20}$ within the range $x\in\left[-2\pi,2\pi\right]$ (see also Appendix B showing how else the Voigt function can be computed).

\begin{figure}[ht]
\begin{center}
\includegraphics[width=24pc]{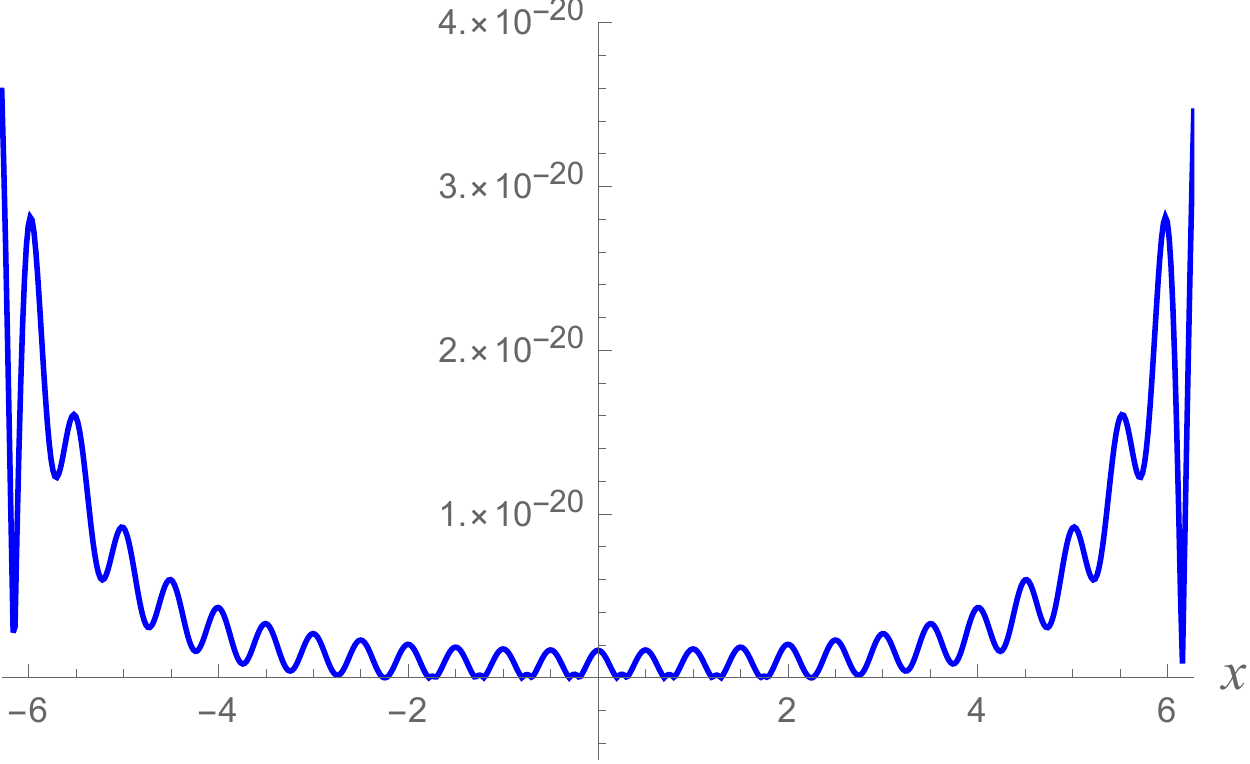}\hspace{2pc}%
\begin{minipage}[b]{28pc}
\vspace{0.3cm}
{\sffamily {\bf{Fig. 7.}} Absolute difference between the original Voigt function $K\left(x,y\right)$ and its approximation \eqref{eq_19} at $y=1$, $a=2$, $M = 6$, $N = 55$, $h=0.078$ and $\sigma = 5$.}
\end{minipage}
\end{center}
\end{figure}

This integration method may be promising in computation of the more sophisticated functions, known as the Voigt function moment integrals, that are used in the atmospheric radiative transfer model MODTRAN \cite{Berk2013}. Our preliminary results show that it can also be efficient in derivation of the rapid and accurate approximation series for some of the integrals involving the Gaussian function $e^{-t^2}$ reported in the recent paper \cite{Chesneau2018}.

\section{Conclusion}

We generalize a methodology shown in our earlier publication \cite{Abrarov2015a} and show as an example how to derive a rational approximation of the sinc function $\text{sinc}\left(\pi \nu \right)$ by sampling and the Fourier transforms. Despite that the sinc function is not easy to approximate, our results reveal that with only $32$ summation terms the absolute difference between the sinc function $\text{sinc}\left(\pi \nu \right)$ and its rational approximation does not exceed $3.2 \times 10^{-3}$ within the range $t \in \left[-2\pi,2\pi\right]$.

\section*{Acknowledgments}

This work is supported by National Research Council Canada, Thoth Technology Inc. and York University. The authors wish to thank the Reviewers for constructive comments and recommendations. Useful remarks and suggestions from John M. Campbell (York University) are greatly appreciated.

\bigskip
\newpage
\section*{Appendix A}

\footnotesize
\begin{verbatim}

function rationapp(opt)

% This function file computes a rational approximation of the sinc and
% v*exp(-(v)^2) functions by using equation (16).

if nargin == 0
    opt = 1; % choose option 1 if argument is missing
    disp('Missing input parameter! Option opt = 1 is assigned.')
end

switch opt
    case 1
        % Define parameters
        a = 0.6;
        k = 35;
        sigma = 2.7;
        M = 6;
        h = 0.04;
        N = 28;
    
        n = 0:N; % define array n
    
        f = 1./((2*(h*n - a)).^(2*k) + 1); % function f(t - a)
        % f = exp(-(2*(n*h - a)).^(2*k)); % alternative function f(t - a)
    
        func = 'sinc(nu)'; % string for first function
        str = 'Sinc approximation'; % string for y-label

    case 2
        % Define parameters
        a = 2;
        sigma = 5;
        M = 6;
        h = 0.078;
        N = 55;

        n = 0:N; % define array n
    
        f = pi^(3/2)*1i*(n*h - a).*exp(-(pi*(n*h - ...
            a)).^2); % function f(t - a)
    
        func = 'nu.*exp(-(nu).^2)'; % string for second function
        str = 'Approximation {\it{\nu e^{-\nu^2}}}';  % string for y-label
        
    otherwise
        disp(['Parameter ',num2str(opt),' is wrong! Choose eather 1 or 2.'])
        return
end

m = exp(sigma*n*h); % exponential multiplier
fm = f.*m;
M2=2^(M-1);

% Define column arrays
alpha = zeros(M2,1);
beta = zeros(M2,1);

% Compute the expansion coefficients
% Equations (13c), (14a) and (14b)
gamma=pi*(2*[1:M2]'-1)./(2^M*h);
for m = 1:M2
alpha(m) = 1/M2*sum(fm.*cos(gamma(m).*n*h),2);
beta(m) = 1/M2*sum(fm.*gamma(m).*sin(gamma(m).*n*h),2);
end

%--------------------------------------------------------------------------
% APPROXIMATION (16)
%--------------------------------------------------------------------------
nu = linspace(-2*pi,2*pi,1000); % define array for the argument nu
funcApp = 0; % initiate function approximation as zero

for m=1:M2
funcApp = funcApp + (alpha(m)*(sigma + 2*pi*1i*nu) + ...
beta(m))./(gamma(m)^2 + (sigma + 2*pi*1i*nu).^2);
end
funcApp = exp(2*pi*1i*nu*a).*funcApp; % final result
%--------------------------------------------------------------------------

% FIGURE 1
figure1 = figure;
axes1 = axes('Parent',figure1,'FontSize',12);
xlim(axes1,[-2*pi,2*pi]);
box(axes1,'on');
grid(axes1,'on');
hold(axes1,'all');
plot1 = plot(nu,[real(funcApp);eval(func)],'Parent',axes1);
set(plot1(1),'LineWidth',3,'Color',[0 1 1]);
set(plot1(2),'LineStyle','--','LineWidth',2);
xlabel('Parameter \it{\nu}','FontSize',14);
ylabel(str,'FontSize',14);

% FIGURE 2
figure2 = figure;
axes2 = axes('Parent',figure2,'FontSize',12);
xlim(axes2,[-2*pi,2*pi]);
box(axes2,'on');
grid(axes2,'on');
hold(axes2,'all');
plot2 = plot(nu,abs(eval(func) - real(funcApp)),'Parent',axes2);
set(plot2(1),'LineWidth',1,'Color',[0 0 0]);
xlabel('Parameter \it{\nu}','FontSize',14);
ylabel('Absolute difference','FontSize',14);

\end{verbatim}
\normalsize

\bigskip
\newpage
\section*{Appendix B}

Due to symmetric forms of the equations \eqref{eq_2} and \eqref{eq_3} an approximation for the inverse Fourier transform follows immediately from the equation \eqref{eq_16} as
\[
\mathcal{F}^{-1}\left\{ F \left( \nu \right) \right\} \left( t \right) \approx e^{-2\pi i t a}\sum\limits_{m = 1}^{{2^{M - 1}}} {\frac{{{\alpha _m^*} \left({\sigma - 2\pi i t} \right) + {\beta _m^*}}}{{\gamma _m^2 + {{\left(\sigma - 2\pi i t \right)}^2}}}},
\]
where the expansion coefficients can be found accordingly as
\[
{\alpha _m^*} = \frac{1}{{{2^{M - 1}}}}\sum\limits_{n = 0}^N { {F\left( {nh - a} \right){e^{\sigma nh}}\cos \left( {{\gamma _m}n h} \right)} }
\]
and
\[
{\beta _m^*} = \frac{1}{{{2^{M - 1}}}}\sum\limits_{n = 0}^N {{F\left( {nh - a} \right){e^{\sigma nh}}{\gamma _m}\sin \left( {{\gamma _m}n h} \right)} }.
\]
Therefore, the Voigt function \eqref{eq_18} can also be computed by using the contour integral of kind
\[
\begin{aligned}
&\hspace{-1cm}\oint\limits_C {\mathcal{F}^{-1}\left\{ {F\left( \nu \right)} \right\}\left(t\right)\, h\left(t,x,y,z, \dots \right) d t } \\
& \hspace{0.5cm} = \oint\limits_C { \left( e^{-2\pi i t a}\sum\limits_{m = 1}^{{2^{M - 1}}} {\frac{{{\alpha _m^*} \left(\sigma - {2\pi i t} \right) + {\beta _m^*}}}{{\gamma _m^2 + {{\left(\sigma - {2\pi i t} \right)}^2}}}} \right)} h\left(t,x,y,z, \dots \right)\,d t,
\end{aligned}
\]
such that the functions are $F(\nu) = \sqrt{\pi } e^{-\pi ^2 \nu ^2} \Leftrightarrow \mathcal{F}^{-1} \left\{\sqrt{\pi } e^{-\pi ^2 \nu ^2} \right\}\left(t\right) = e^{-{t^2}}$ and 
$$
h\left(t,x,y\right) = \frac{y}{\pi \left( y^2 + \left( x - t \right)^2 \right)}.
$$


\end{document}